\documentclass[11pt]{article}
\newcommand{\subsetm}{\stackrel{{\rm m}}{\subset}}

\newcommand{\Z}{\mathbb Z}
\newcommand{\R}{\mathbb R}

\newtheorem{theorem}{Theorem}[section]

\newtheorem{lemma}[theorem]{Lemma}

\newtheorem{definition}[theorem]{Definition}

\newtheorem{example}[theorem]{Example}

\newtheorem{remark}[theorem]{Remark}

\usepackage{amsfonts,amssymb}

\input{epsf}

\date{}

\begin{document}

\title{ Generalizations of Quandle Cocycle Invariants and Alexander Modules 
from Quandle Modules}

\author{
J. Scott Carter\footnote{Supported in part by NSF Grant DMS \#0301095.}
\\University of South Alabama \\
Mobile, AL 36688, U.S.A. \\ carter@jaguar1.usouthal.edu
\and 
Masahico Saito\footnote{Supported in part by NSF Grant DMS \#0301089.}
\\ University of South Florida
\\ Tampa, FL 33620,  U.S.A.  \\ saito@math.usf.edu
}

\maketitle

\section{Introduction} \label{intro}

Quandle cohomology theory 
was developed \cite{CJKLS} to define invariants
of classical knots and knotted surfaces in state-sum form,
called  quandle cocycle (knot) invariants.
The quandle cohomology theory is a modification
of rack cohomology theory  which was defined in \cite{FRS}.
The cocycle knot invariants are analogous in their definitions
to the Dijkgraaf-Witten invariants
\cite{DW} of triangulated $3$-manifolds with finite gauge groups, 
but they 
use quandle knot colorings as spins and cocycles as 
Boltzmann weights.
 In \cite{CEGS},  the quandle cocycle invariants were generalized
in three different directions, using generalizations of 
quandle homology theory provided by  Andruskiewitsch and Gra\~{n}a  
\cite{AG}, which is compared to the group cohomology theories with 
the group actions on the coefficient groups.
This paper is a written version of our talk given at 
Intelligence of Low Dimensional Topology in Shodo-Shima.
It is a short summary of \cite{CEGS} with 
some results from \cite{ribbon} and a few new observations.
We would like to thank the organizers for holding such an exciting conference in
a beautiful location.

\section{Preliminary: Quandles and colorings} \label{prelimsec}

A {\it quandle}, $X$, is a set with a binary operation 
$(a, b) \mapsto a * b$
such that

(I) For any $a \in X$,
$a* a =a$.

(II) For any $a,b \in X$, there is a unique $c \in X$ such that 
$a= c*b$.

(III) 
For any $a,b,c \in X$, we have
$ (a*b)*c=(a*c)*(b*c). $

\noindent
A {\it rack} is a set with a binary operation that satisfies 
(II) and (III).
Racks and quandles have been studied in, for example, 
\cite{AG,Br88,FR,Joyce,K&P,Matveev}.


The following are typical examples of quandles.  
A group $X=G$ with
conjugation
as the quandle operation:
$a*b=b a b^{-1}$. 
We denote by Conj$(G)$ the quandle defined for a group $G$ by 
$a*b=bab^{-1}$. 
Any subset of $G$ that is closed under such conjugation 
is also a quandle. 

Any $\Lambda (={\Z }[t, t^{-1}])$-module $M$  
is a quandle with
$a*b=ta+(1-t)b$, $a,b \in M$, 
that is 
called an {\it  Alexander  quandle}.
Let $n$ be a positive integer, and 
for elements  $i, j \in \{ 0, 1, \ldots , n-1 \}$, define
$i\ast j \equiv 2j-i \pmod{n}$.
Then $\ast$ defines a quandle
structure  called the {\it dihedral quandle},
  $R_n$.
This set can be identified with  the
set of reflections of a regular $n$-gon
  with conjugation
as the quandle operation; it also is isomorphic to an Alexander quandle 
${\Z }_n[t, t^{-1}]/(t+1)$. 
As a set of reflections of the 
regular 
$n$-gon, 
 $R_n$ can be considered as a subquandle of ${\mbox{\rm Conj}}(\Sigma_n)$.

\begin{figure}
\begin{center}
\mbox{
\epsfxsize=2.5in
\epsfbox{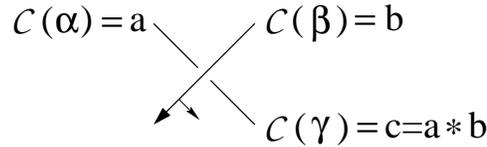} 
}
\end{center}
\caption{ Quandle relation at a crossing  }
\label{qcolor} 
\end{figure}

Let $X$ be a fixed quandle.
Let $K$ be a given oriented classical knot or link diagram,
and let ${\cal R}$ be the set of (over-)arcs. 
The normals are given in such a way that (tangent, 
normal) agrees with
the orientation of the plane, see Fig.~\ref{qcolor}. 
A (quandle) {\it coloring} ${\cal C}$ is a map 
${\cal C} : {\cal R} \rightarrow X$ such that at every crossing,
the relation depicted in Fig.~\ref{qcolor} holds. 
More specifically, let $\beta$ be the over-arc at a crossing,
and $\alpha$, $\gamma$ be under-arcs such that the normal of the over-arc
points from $\alpha$ to $\gamma$.
(In this case, $\alpha$ is called the {\it source arc} and  $\gamma$ 
is called the {\it target arc}.) 
Then it is required that ${\cal C}(\gamma)={\cal C}(\alpha)*{\cal C}(\beta)$.
The  colors ${\cal C}(\alpha)$, ${\cal C}(\beta)$
are called {\it source} colors.

\section{Quandle Modules}

We recall some information from \cite{AG}, 
but with notation changed to match our conventions. 

Let $X$ be a quandle. Let $\Omega (X)$ be the free ${\Z}$-algebra 
generated by 
$\eta_{x,y},$ $\tau_{x,y}$   for $x,y \in X$ such that 
$\eta_{x,y}$ is invertible for every $x,y \in X$.  
Define $\Z (X)$ to be the  quotient $\Z(X)= \Omega(X)/ R$ 
where $R$ is the ideal 
generated by 

\begin{enumerate}
\setlength{\itemsep}{-2pt}
\item \hfil $ \eta_{x*y,z}\eta_{x,y} -  \eta_{x*z,y*z}\eta_{x,z}$ \hfill \hfill
\item \hfil $ \eta_{x*y,z}\tau_{x,y} -  \tau_{x*z,y*z}\eta_{y,z}$ \hfill \hfill
\item \hfil $\tau_{x*y,z}- \eta_{x*z,y*z} \tau_{x,z}- 
 \tau_{x*z,y*z}\tau_{y,z}$ \hfill \hfill
\item \hfil $\tau_{x,x} + \eta_{x,x} -1$ \hfill \hfill
\end{enumerate}

The algebra  $\Z(X)$ thus defined is called the  {\it quandle algebra} 
over $X$. 
In $\Z(X)$, we define elements 
$\overline{\eta_{z,y}} =\eta^{-1}_{z\overline{*}y,y}$
and 
$\overline{\tau_{z,y}} = -\overline{\eta_{z,y}}\tau_{z\overline{*}y,y}$.

A {\it representation} of $\Z(X)$ is an an algebra homomorphism  
$\Z(X) \rightarrow {\mbox{\rm End}} (G)$, and we denote the 
image of the generators by the same symbols. Given a representation of 
$\Z(X)$ we say that $G$ is a $\Z(X)$-module, or {\it a quandle module}. 
The action of $\Z(X)$ on $G$ is written by the left action,  
and denoted by  $(\rho, g) \mapsto \rho g  (= \rho \cdot g = \rho (g) )$,
for  $g\in G$ and  $\rho \in End(G).$

\begin{figure}[h]
\begin{center}
\mbox{
\epsfxsize=2.5in
\epsfbox{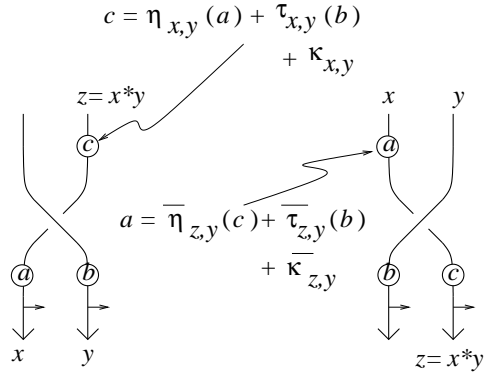} 
}
\end{center}
\caption{The geometric notation at a crossing}
\label{gencolor} 
\end{figure}


Diagrammatic conventions of $\eta$ and $\tau$ are depicted in 
Fig.~\ref{gencolor}, where $\kappa$ is a generalized $2$-cocycle that will
appear later in this paper. 
We have to leave details to \cite{CEGS}.

\begin{example} {\bf \cite{AG}} \label{AGexample} {\rm 
Let $\Lambda = \Z[t,t^{-1}]$ denote the ring of Laurent polynomials. 
Then any  $\Lambda$-module $M$ is a $\Z(X)$-module for any quandle $X$,
by $\eta_{x,y} (a)=ta$ and $\tau_{x,y} (b) = (1-t) (b) $
for any $x,y \in X$.

The group  $G_X=\langle x \in X \ | \ x*y=yxy^{-1} \rangle$
is called the {\it enveloping group} \cite{AG}
(and the {\it associated group} in \cite{FR}). 
For any quandle $X$,
any $G_X$-module $M$  is a  $\Z(X)$-module by 
 $\eta_{x,y} (a)=y a$ and $\tau_{x,y} (b) = (1- x*y) (b) $,
where $x, y \in X$, $a, b \in M$.
} \end{example} 

\begin{remark} {\rm 
The second example above is the Fox's free
derivative of the  braid group representation into 
the automorphism group of the free group which corresponds to
 $(x, y)\mapsto (y, yxy^{-1})$. 
Indeed, one computes 
$$ \eta_{x,y} = \frac{\partial}{\partial x} (yxy^{-1})=y,
\quad \tau_{x,y} = \frac{\partial}{\partial y} (yxy^{-1})= 1- yxy^{-1}. $$

Wada 
listed certain types of braid group representations.
Among them are 
of the form that a standard braid generator 
acts as $(x,y) \mapsto (y, w(x,y))$, where 
$w(x,y)=y^m x y^{-m}$ for some integer $m$, 
and $w(x,y)=yx^{-1}y$. 
We remark here that the Fox's derivative of these 
give rise to quandle module structures as well.
Cocycle invariants for these quandle module structures
would be of future research interest.
} \end{remark}

\section{Generalized quandle homology theory}
Consider 
the free 
right 
$\Z ( X )$-module $C_n(X)=  \Z ( X ) X^n$ 
with basis $X^n$ (for $n=0$, $X^0$ is a singleton $\{ x_0 \}$, 
for a fixed element $x_0 \in X$). 
In \cite{AG},  boundary operators
 $\partial =\partial_n: C_{n+1}(X) \rightarrow C_n(X)$ 
are defined by 
\begin{eqnarray*} \lefteqn{\partial (x_1, \ldots, x_{n+1} ) } \\
 &=& (-1)^{n+1} {\displaystyle \sum_{i=2}^{n+1} (-1)^i 
\eta_{ [x_1, \ldots, \widehat{x_i}, \ldots, x_{n+1} ], [x_i, \ldots, x_{n+1}] }
(x_1, \ldots, \widehat{x_i}, \ldots, x_{n+1} ) } \\
 & & \\
 & & - {\displaystyle (-1)^{n+1} \sum_{i=2}^{n+1} (-1)^i 
(x_1*x_i, \ldots, x_{i-1}*x_i, x_{i+1},  \ldots, x_{n+1} ) }  \\
 & & \\
 & & + (-1)^{n+1} 
 \tau_{[x_1, x_3, \ldots, x_{n+1}], [x_2, x_3, \ldots, x_{n+1}]}
(x_2, \ldots, x_{n+1} ) , \\[3mm]
\mbox{where} & & 
 [x_1, x_2, \ldots, x_n] = ( ( \cdots ( x_1 * x_2 ) * 
x_3  ) * \cdots ) * x_n
\end{eqnarray*}
for $n > 0$, and 
$\partial_1(x) = - \tau_{x \bar * x_0, x_0}$ for $n=0$.
The notational conventions  are slightly different from \cite{AG}.   
In particular, the  $2$-cocycle 
condition for a $2$-cochain $\kappa_{x,y}$ 
in this homology theory 
is written as 
$$ \eta_{x*y, z} (\kappa_{x,y}) + \kappa_{x*y, z}
= \eta_{x*z, y*z} (\kappa_{x,z})  +  \tau_{x*z, y*z} (\kappa_{y,z})
 + \kappa_{x*z,y*z}, $$ 
for any $x,y,z \in X$.
We call this a {\it generalized $($rack$)$ $2$-cocycle condition}.
When $\kappa $ further satisfies $\kappa_{x,x}=0$ for any $x \in X$,
we call it a {\it generalized quandle $2$-cocycle.}

\section{Assigning homology classes to colored diagrams}  
Here we review only the knotted surface case.
The classical case is similar and the triple points are 
replaced by crossings and 
$3$-cocycles are replaced by $2$-cocycles.
This method was independently developed in \cite{Tanaka}. 

\begin{figure}[htb]
\begin{center}
\mbox{
\epsfxsize=3.5in
\epsfbox{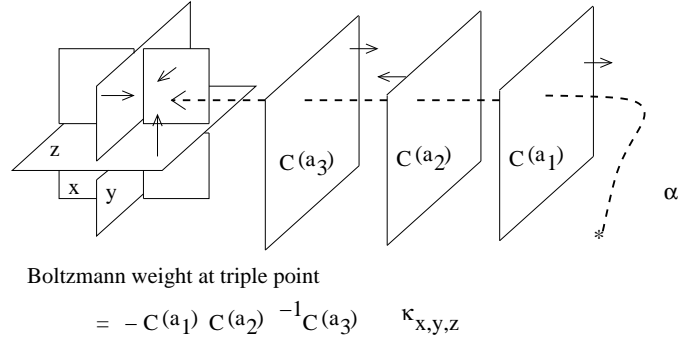} 
}
\end{center}
\caption{The weight at a triple point}
\label{trptweight} 
\end{figure}

A  diagram $D$ of a knotted surface $K$ is given in  $3$-space. 
We assume the surface is oriented and 
use orientation normals to indicate the orientation. 
In a neighborhood of each triple point, there are eight regions 
that are separated by the sheets of the surface since
 the triple point looks like the intersection of the 
$3$-coordinate planes in some parametrization.
The region into which 
all normals point
is called the {\it target} region.  
Let $\gamma$ be an arc from the region at infinity of  $3$-space
 to the target region of a given triple point $r$. Assume  that 
$\gamma$
intersects
$D$ transversely in a finitely many points thereby missing 
double point curves, branch points, and triple points. 
Let $a_i$, $i=1, \ldots, k$, in this order, be the sheets of $D$
that intersect $\gamma$ from the region at infinity to the triple point $r$.
Let ${\cal C}$ be a coloring of $D$ by a fixed finite quandle $X$.
See Fig.~\ref{trptweight}.

\begin{definition} {\rm 
The $3$-chain 
$$ {\cal C}(D) = \sum_r \epsilon(r)  
({\cal C}(a_1)^{\epsilon (a_1) } {\cal C}( a_2)^{\epsilon (a_2) }
 \cdots{\cal C}( a_k)^{\epsilon (a_k) } ) 
 (x, y, z) \in C_3(X; \Z G_X) $$
is called the {\it $3$-chain represented by the diagram $D$ with the coloring 
${\cal C}$}.

} \end{definition}

The following  is proved from definitions.

\begin{lemma} 
For any knotted surface diagram $D$ with a coloring ${\cal C}$, 
the $3$-chain $ {\cal C}(D)$
  represented by the diagram $D$ with the coloring 
${\cal C}$ is a $3$-cycle: $ {\cal C}(D) \in Z_3(X; \Z G_X) $.
\end{lemma}



\section{Cocycle invariants }

We continue with the surface case, as the classical case is similar.

Let $X$ be a finite quandle and $\Z (X)$ be its quandle algebra
with generators $\{ \eta _{x,y}^{\pm 1}  \}_{ x, y \in X} $
and $\{ \tau _{x,y}  \}_{ x, y \in X} $.
Let $G$ be 
an abelian group that is a $G_X$-module. Recall that  
this induces a 
$\Z (X)$-module structure  
given by $\eta_{x,y} g = y g$ and $\tau_{x,y}(g) = (1 - x*y)g$
 for $g \in G$ and 
$x,y \in X$. 
Let $\kappa_{x,y,z}$ be a generalized quandle $3$-cocycle of $X$ with 
the coefficient group $G$. 
Thus the generalized $3$-cocycle condition, in this setting, is 
written as 
\begin{eqnarray*}
\lefteqn{
 w \kappa_{x,y,z} + \kappa_{x*z, y*z, w} + ((y*z)*w) \kappa_{x,z,w}
+ \kappa_{y,z,w} } \\
&= & (((x*y)*z))*w)  
\kappa_{y,z,w}  + \kappa_{x*y, z, w} 
+ (z*w) \kappa_{x,y,w} + \kappa_{x*w, y*w, z*w}.
\end{eqnarray*}
We require further that 
$\kappa_{x,x,y}=\kappa_{x,y,y}=0$.
These conditions are called {\it quandle cocycle conditions},
and a $3$-cocycle that satisfies the quandle cocycle conditions
is called a (generalized) {\it quandle $3$-cocycle}.
A cocycle invariant of knotted surfaces
 will be defined using such a $3$-cocycle. 



\begin{definition}
{\rm The {\em Boltzmann weight} $B( {\cal C}, r,\gamma)$
for the triple point $r$, for a coloring ${\cal C}$, with respect to $\gamma$,
is defined by 
$$ B( {\cal C}, r,\gamma)=
 \epsilon(r)  ({\cal C}(a_1)^{\epsilon (a_1) } {\cal C}( a_2)^{\epsilon (a_2) }
 \cdots{\cal C}( a_k)^{\epsilon (a_k) } ) 
 \kappa_{x, y, z} \in  G, $$
where $x,y,z$ are the  {\it color triplet} at the given triple point $r$
($x$ is assigned on the bottom sheet from which the normals of the
middle and top sheets
point, and $y$ is assigned to the middle sheet from which the normal
of the top sheet points, and $z$ is assigned to the top sheet). 
The sign $ \epsilon (r)$ 
is the sign of the triple point $r$.
The exponent $\epsilon (a_j)$ is $1$ is the arc $\gamma$
crosses the arc $a_j$ against its normal, and is $-1$ otherwise,
for $j=1, \ldots, k$. 
} \end{definition}


\begin{definition} {\bf  \cite{CEGS}} {\rm
The family 
$\Phi_{\kappa} (K) = \{ \sum_{r} B( {\cal C}, r) \}_{{\cal C} \in \mbox{\rm Col}_{\rm X} (D) } $
is called the quandle cocycle invariant with respect to 
the (generalized) $3$-cocycle $\kappa$. 
} \end{definition}

\begin{theorem} {\bf \cite{CEGS}}
The family 
$\Phi_{\kappa} (K) $ does not depend on the choice of a diagram $D$
of a given knotted surface $K$, so that it is a well-defined knot invariant.
\end{theorem}

The cocycle invariant can be regarded as a family over all colorings 
(or the formal sum) of the Kronecker 
product $B( {\cal C}, r) =\langle \kappa, {\cal C}(D) \rangle$.

\section{Computations}

The generalized cocycle invariants were computed using
{\it Maple} and {\it Mathematica} in \cite{CEGS}. 
Here we include the table of the invariant for a certain $3$-cocycle
of $R_3$ with coefficient group $\Z^3$ where 
$R_3$ acts as permutations of factors of the vectors in  $\Z^3$.
The table is for the $2$-twist spin of classical knots in the table, 
up to 8 crossings. Those for up to 9 crossing were computed in \cite{CEGS}.
The $3$-cocycle used has two free variables $q_1$ and $q_2$, 
so that the values in the table contain these. 
The notation $\sqcup_n$ indicates $n$ copies of the vector.

\begin{table}[h] 
\begin{center}
{\begin{tabular}{||l||l||}\hline \hline
Knot $K$ & $\Phi_{\kappa}({\rm Tw}^2 (K) )$
\\ \hline \hline 
$3_1$ & $\sqcup_9 (0,0,0). $ 
\\ \hline 
$6_1$ & $\sqcup_3 (0,0,0),\; 
(-q_2, 0, q_2),\;  (q_2, q_1,  - q_1-q_2),\; 
     (q_1, 0, -q_1),\; $
\\ \hline 
 & $  (0, -q_1, q_1),\; 
      ( - q_1+q_2,  q_1 -q_2, 0),\;  (-q_2,  - q_1+q_2, q_1).$  
\\ \hline 
$7_4$  & $ \sqcup_3 (0,0,0),\; 
(-q_2, -2 q_1, 2 q_1 + q_2),\; 
  (-q_1 + q_2, q_1, -q_2),\; $
\\ \hline 
 & $ (-q_1, 2 q_1, -q_1)\; ,
 (-q_1, -q_1, 2 q_1),\;  (q_2, -q_2, 0),\; 
     (-q_2, q_2, 0) . $ 
\\ \hline 
$7_7$    &  $ \sqcup_3 (0,0,0),\;
\sqcup_2 (q_1, 0, -q_1),\; \sqcup_2 (0, -q_1, q_1),\; 
\sqcup_2  (-q_1, q_1, 0).  $ 
\\ \hline 
$8_5$  &  $ \sqcup_9 (0,0,0). $ 
\\ \hline 
$8_{10}$ &  $ \sqcup_9 (0,0,0). $ 
\\ \hline
$8_{11}$   &  $ \sqcup_3 (0,0,0),\; 
 (q_2, 0, -q_2),\;  (-q_2, -q_1,   q_1+q_2),\; 
     (-q_1, 0, q_1),\; $
\\ \hline 
 & $ (0, q_1, -q_1),\; 
  (q_1 - q_2,  - q_1+q_2, 0),\;  (q_2, q_1 - q_2, -q_1). $  
\\ \hline 
$8_{15}$   &  $ \sqcup_3 (0,0,0),\;  
 \sqcup_3 (0, -q_1, q_1),\; \sqcup_3 (-q_1, q_1, 0). $  
\\ \hline 
$8_{18}$    &  $ \sqcup_9 (0,0,0),\;
 \sqcup_6 (q_1, 0, -q_1),\;
 \sqcup_6  (-q_1, q_1, 0),\;
 \sqcup_6 (0, -q_1, q_1).$
\\ \hline
$8_{19}$    &   $ \sqcup_9 (0,0,0). $ 
\\ \hline 
$8_{20} $   &  $ \sqcup_9 (0,0,0).  $  
\\ \hline 
$8_{21}$    &   $ \sqcup_9 (0,0,0). $  
\\ \hline 
\end{tabular} } \end{center}
\caption{A table of cocycle invariants for twist spun knots}
\label{twistspintable1}
\end{table}

\section{Applications}

The following topological applications have been found.

\bigskip

\noindent
{\bf Non-invertibility of knotted surfaces.} 
A knot(ed surface) is called invertible if it is equivalent to itself with the opposite orientation, with the orientation of the space fixed.
 
Fox \cite{Fox61a} presented a  non-invertible knotted sphere
using  asymmetric knot   modules.
Farber \cite{Far75} 
showed that the $2$-twist spun trefoil  was 
non-invertible  using the Farber-Levine pairing
(see also Hillman \cite{Hill81}).   
Ruberman \cite{Rub83}  
used Casson-Gordon invariants  
to prove the same result,
with other new examples of non-invertible knotted spheres. 
Neither technique applies
directly to the same knot with 
trivial $1$-handles attached
(in this case the knot is a surface with a higher genus).
Kawauchi~\cite{Kawa86a,Kawa90a}  
has generalized the Farber-Levine pairing to higher genus surfaces,
showing that such a  surface 
is also non-invertible. 
Gordon~\cite{Gor} showed that a large family of knotted spheres
are indeed non-invertible using fibrations.

In \cite{CJKLS} the original quandle cocycle invariants were used
to detect non-invertibility of the $2$-twist spun trefoil.
Furthermore, 
using cocycle invariants for proving non-invertibility
applies to stabilzed  surfaces.
These are obtained 
by attaching trivial $1$-handles.
Satoh~\cite{Sat01b*} applied this method to prove non-invertibility of 
an infinite family of twist spins of torus knots and their stabilizations.
A similar result has been obtained for $2$-bridge knots by Iwakiri\cite{Iwakiri}. 
{}From the generalized invariant, we have:

\begin{theorem} {\bf \cite{CEGS}}
For any positive integer $k$, 
the $2k$-twist spun  of all the 
$3$-colorable  knots in the table up to $9$ crossings  excluding $8_{20}$, 
as well as their stabilized surfaces of any genus, 
are non-invertible.
\end{theorem}
Some of the computations that yield this result are presented in Table 1.

\bigskip

\noindent
 {\bf Minimal triple point numbers of projections of knotted surfaces.}
Classical knot tables have played a pivotal role in the
history of classical knot theory. 
Knot tables are organized according to (minimal) crossing numbers.
An analogue is the minimal number of triple points among all gereric projections of a knotted surface, called the {\it triple point number}  ${\cal T}(F)$ of a knotted surface $F$.

Progress has been made in \cite{Sat00a, Sat01a*, Shima01*}
about triple point numbers, but 
 there were no examples of knotted spheres whose
triple point numbers 
were concretely determined
until a breakthrough was given in \cite{SatShi01a*,SatShi01b*},
in which the following were proved:
\begin{eqnarray*}
{\cal T}\ ( \ \mbox{$2$-twist spun trefoil} \ ) & = & 4, \\ 
{\cal T}\ (\  \mbox{$3$-twist spun trefoil}\ ) & = & 6.  
\end{eqnarray*}

Since then further applications of quandle colorings 
and the cocycle invariants to the triple point numbers 
have been obtained (\cite{Hatake}, for example).

\begin{remark}{\rm 
These results can be interpreted as  a pseudo-norm on
quandle homology. Specifically, the minimum number of 
generators with which a non-trivial homology class is represented 
is regarded as a pseudo-norm
and gives a lower bound for the triple point numbers.
} \end{remark}

Generalized cocycle invariants often give higher lower bounds.
For example, it can be used to reprove:  

\begin{theorem} {\rm ({\bf  \cite{Kama92b}}, see also {\bf \cite{Kawa}})}%
For any positive integer $n$, there exists a knotted surface $K$
such that ${\cal T}(K) > n$.
\end{theorem}

In fact, the theorem is verified by specific examples.
With $R_3$ and a certain cocycle, it is proved that
if a knotted surface $K$ has a non-trivial 
cocycle invariant, then 
 for any positive integer $n$, there is a positive integer $k$
such that ${\cal T}(\tau^{2k} (K)) \geq n$. 

\bigskip

\begin{figure}[htb]
\begin{center}
\mbox{
\epsfxsize=2in
\epsfbox{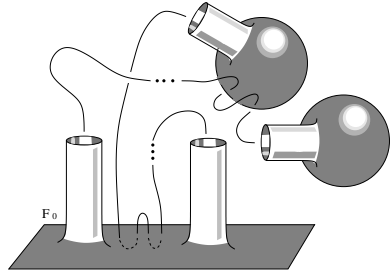} 
}
\end{center}
\caption{Ribbon concordance}
\label{ribbonidea} 
\end{figure}

\noindent
 {\bf Ribbon concordance of knotted surfaces.}
Let $F_0$ and $F_1$ be connected knotted surfaces of the same genus. 
We say that {\it $F_1$ is ribbon concordant to $F_0$} 
if there is a concordance $C$ 
in ${\R}^4\times [0,1]$ between 
$F_1\subset{\R}^4\times\{1\}$ and 
$F_0\subset{\R}^4\times\{0\}$ 
such that the restriction to $C$ 
of the projection ${\R}^4\times[0,1]\rightarrow[0,1]$ 
is a Morse function with critical points 
of index $0$ and $1$ only. 
We write $F_1\geq F_0$. 
Note that if $F_1\geq F_0$, 
then there is a set of $n$ $1$-handles 
on a split union of $F_0$ and $n$ trivial sphere, 
for some $n\geq 0$, 
such that $F_1$ is obtained by surgeries 
along these handles (Fig.~\ref{ribbonidea}).

The notion of ribbon concordance was 
originally introduced by Gordon \cite{Gor} 
for classical knots in ${\R}^3$, 
and there are several studies 
found in \cite{Gil, Miy1, Miy2, Sil}, 
for example. 

Given knotted surfaces $F_0$ and $F_1$, 
it is natural to ask whether 
$F_1$ is ribbon concordant to $F_0$. 
Cochran \cite{Coc} gave 
a necessary condition 
for a sphere-knot $F$ to be ribbon 
in terms the knot group $\pi_1({\R}^4\setminus F)$. 
In \cite{ribbon}, new necessary conditions  were given for 
a pair of  knotted surfaces  to be ribbon concordant 
by using quandle cocycle invariants.

The cocycle invariant $\Phi_{\theta}(F)$ is 
regarded as 
a multi-set of elements in the coefficient group $A$ 
of the cohomology 
where repetitions of the same element are allowed. 
For two multi-sets $A'$ and $A''$ of $A$, 
we use the notation $A'\subsetm A''$ 
if for any $a\in A'$ it holds that $a\in A''$. 
In other words, 
$A'\subsetm A''$ if and only if 
$\tilde{A'}\subset \tilde{A''}$ 
where $\tilde{A'}$ and $\tilde{A''}$ are 
the subsets of $A$ obtained from $A'$ and $A''$ 
by eliminating the multiplicity of elements, 
respectively. 

\begin{theorem} {\bf \cite{ribbon} } \label{thm11} 
If $F_1\geq F_0$, 
then $\Phi_{\theta}(F_1)\subsetm\Phi_{\theta}(F_0)$. 
\end{theorem} 

By Theorem~\ref{thm11}, 
we give many examples of pairs of  knotted surfaces 
such that one is {\it not} 
ribbon concordant to another. 

To generalize this result to surfaces without triple points, 
a new cocycle invariants defined on
$H_1(F)$ for a surface $F$ was constructed in \cite{ribbon} as well.

\begin{remark} 
{\rm Kawauchi points out that 
the linking signature of a certain family of surfaces 
is invariant under ribbon concordance. 
This result has not appeared in any paper, 
but can be obtained as a corollary of \cite{Kaw}. 
}
\end{remark}

\section{Module invariants and twisted Alexander invariants}

A braid word $w$ (of $k$-strings), or a $k$-braid word,  is 
a product of standard  generators $\sigma_1, \ldots, \sigma_{k-1}$ 
of the braid group ${\cal B}_k$ of $k$-strings and their inverses.
A braid word $w$ represents an element $[w]$ of the braid group 
${\cal B}_k$.
Geometrically, $w$ is represented by a diagram in a rectangular box
with  $k$ end points  at the top, and $k$ end points at the bottom,
where the strings   go down monotonically. 
Each generator or its inverse is represented by a crossing in a diagram.
We use the same letter $w$ for a choice of such a diagram. 
Let $\hat{w}$  denote the closure of 
the 
diagram 
$w$. 
Quandle colorings of 
$w$ 
are  defined in 
exactly the same manner as in the case of knots.
However, the quandle elements at the top and the bottom of a diagram of 
$w$ do not necessarily coincide. 
But when when we consider a coloring of the  link  $\hat{w}$,
the 
 quandle elements at the top and the bottom of a diagram of 
$w$ do coincide.

\begin{figure}[h]
\begin{center}
\mbox{
\epsfxsize=1in
\epsfbox{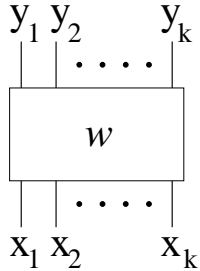} 
}
\end{center}
\caption{A quandle coloring of a braid word $w$}
\label{bcolors} 
\end{figure}

Let $X$ be a quandle.
Let $\gamma_1, \ldots, \gamma_k$ be the bottom arcs of $w$.
For a given vector $\vec{x}=(x_1, \ldots, x_k) \in X^k$,
assign these elements $ x_1, \ldots, x_k$ on $\gamma_1, \ldots, \gamma_k$ 
 as their colors, respectively. Then from the definition, 
a coloring ${\cal C}$ of $w$ by $X$ is uniquely determined
such that  
${\cal C}(\gamma_i)=x_i$, $i=1, \ldots, k$.
We call such a coloring ${\cal C} $ the {\it  coloring induced from} $\vec{x}$.
Let $\delta_1, \ldots, \delta_k$ be the arcs at the top.
Let 
$\vec{y}=(y_1, \ldots, y_k)=({\cal C}(\delta_1), \ldots, {\cal C}(\delta_k) )\in X^k$ 
be the colors assigned to the top arcs. 
Denote
 this situation by a left action, 
$\vec{y} =  w \cdot \vec{x}$. The colors $\vec{x}$ and $\vec{y}$ are 
called {\it bottom} and {\it top colors}
or {\it color vectors}, 
respectively.
See Fig.~\ref{bcolors}. 

Let $X$ be a quandle and $G$ be a quandle module. 
For 
$\alpha= \eta + \tau $ which acts on $(a,b) \in G^2$ by
$ \alpha_{x,y}(a,b) =\eta_{x,y} (a )  +  \tau_{x,y}(b)$
for any $(x,y) \in X^2$.
Let 
$\tilde{X}=G \times_{\alpha} X$ be the quandle defined by
$(a, x)*(b,y)=(\alpha_{x,y}(a,b), x*y)$ (which is called the 
dynamical extension). 
If $\vec{r} =( (a_1, x_1), \ldots, (a_k, x_k)) $
and $\vec{s} =( (b_1, y_1), \ldots, (b_k, y_k)) \in \tilde{X}^k$
are bottom and top colors of $w \in {\cal B}_k$ by $\tilde{X}$,
respectively, then we write this situation by 
$\vec{b} =  M(w, \vec{x}) \cdot \vec{a}$, 
where $\vec{a}=(a_1, \ldots, a_k)$, 
$\vec{b}=(b_1, \ldots, b_k) \in G^k$. 
Thus $M(w, \vec{x})$ represents a map 
$M(w, \vec{x}): G^k \rightarrow G^k$.

\begin{lemma} {\bf \cite{CEGS}} \label{indeplemma}
If $[w]=[w'] \in {\cal B}_k$, then 
$M(w, \vec{x})=M(w', \vec{x}) : G^k \rightarrow G^k$. 
\end{lemma}

 We call the map 
$M(-, \vec{x}): {\cal B}_k \rightarrow \mbox{Map}(G^k, G^k)$ 
a {\it colored representation}.

\begin{theorem} {\bf \cite{CEGS} }
Let $L$ be a link represented as a closed braid $\hat{w}$,
where $w$ is a $k$-braid word,   and 
${ \mbox{Col}_X(L) }$ be the set of colorings of $L$ by a quandle $X$.
For ${\cal C} \in  \mbox{Col}_X(L)$, 
let $\vec{x}$ be the color vector of bottom strings of 
$w$ that is the restriction of ${\cal C}$. 
Then the family  
$$
{\cal M}(X, \alpha\ ; L)=\{ 
 G^k 
/\mbox{\rm Im} (M(w,\vec{x})-I)  \}_{{\cal C} \in  {\rm Col}_X(L)}$$
of isomorphism classes of modules presented by 
the maps $(M(w,\vec{x})-I)$, 
where $I$ denotes the identity,  is independent of choice of $w$
that represents $L$ as its closed braid,
and thus defines a link invariant.
\end{theorem}

This theorem implies that the following is well-defined.

\begin{definition} {\bf \cite{CEGS}} {\rm 
The family of modules
${\cal M}(X, \alpha \ ; L)
=\{ G^k 
/\mbox{\rm Im}(M(w,\vec{x})-I)\}_{{\cal C} \in {\rm Col}_X(L) } $
is called the {\em quandle module invariant}.
} \end{definition}

This invariant is  related to the twisted Alexander invariant \cite{KL,Lin,Wada}.
Let 
 $ \pi = \pi_1 (S^3 \setminus K  ) =\langle x_1, \ldots, x_s | r_1, \ldots, r_k \rangle$  be a 
 Wirtinger presentation of a classicak knot  $K$, 
 so that relations are  of the form $x_\ell=x_s x_j x_s^{-1}$.
 Let 
  $V=\Z^n $  and 
    $\rho : \pi_1  \rightarrow  GL(V) $  be
     a representation, and 
  $ \epsilon: \pi_1  \rightarrow \Z $
  be the abelianization.
   Finally, let
   $\Z[F_s]$ be    the free group ring generated by $x_1, \ldots, x_s$
    and  
   $\chi : \Z[F_s] \rightarrow \Z[\pi] \rightarrow M_n(\Z [t, t^{-1}])$ be
    the map that  is determined by 
    $\rho \otimes \epsilon : \Z[\pi] \rightarrow M_n(\Z [t, t^{-1}])$
   that sends $\gamma$ to $t^{{\epsilon}(\gamma) }\rho(\gamma)$.
     Then the  module with the presentation $(sn \times kn)$-matrix
$\displaystyle  \left[  \chi \left( \frac{\partial r_i}{\partial x_j } \right) \right] $
is used to define the twisted Alexander invariant.
  

On the other hand, we can define 
a quandle module structure by using 
 the action of $G_X$ on $(\Z[t, t^{-1}])^n$
 defined by  $\rho \otimes \epsilon : \Z[\pi] \rightarrow M_n(\Z [t, t^{-1}])$.
 Then the quandle module invariant has a similar presentation matrix as above.
 Therefore a detailed analysis  of the roles played by the matrix
 $\displaystyle  \left[  \chi \left( \frac{\partial r_i}{\partial x_j } \right) \right] $
will yield  connections among these subjects.

 

\end{document}